\documentclass[12pt,a4paper]{amsart}
\usepackage[utf8]{inputenc}	
\usepackage[T1]{fontenc}
\usepackage[in,headings]{fullpage}
\usepackage{amsbsy, amscd, amsfonts, amsmath, amsrefs, amssymb, amsthm}
\usepackage{graphicx}
\usepackage{float}
\usepackage{color}   
\usepackage[colorlinks=true,linkcolor=blue,citecolor=blue,urlcolor=blue]{hyperref}
\usepackage{comment}
\excludecomment{A}

\newtheorem{theorem}{Theorem}

\newtheorem{proposition}[theorem]{Proposition}

\theoremstyle{definition}

\newtheorem{remark}[theorem]{Remark}

\theoremstyle{remark}

\newcommand{\C}{\mathbf{C}}

\renewcommand{\Re}{\mathop{\mathrm{Re}}\nolimits}
\renewcommand{\Im}{\mathop{\mathrm{Im}}\nolimits}
\newcommand{\Rzeta}{\mathop{\mathcal R }\nolimits}

\newfont{\cmbsy}{cmbsy10}
\newfont{\cmmib}{cmmib10}
\newcommand{\Orden}{\mathop{\hbox{\cmbsy O}}\nolimits}
\newcommand{\orden}{\mathop{\hbox{\cmmib o}}\nolimits}

\begin{document}

\title{Note on the asymptotic of the auxiliary function}
\author[Arias de Reyna]{J. Arias de Reyna}
\address{%
Universidad de Sevilla \\ 
Facultad de Matem\'aticas \\ 
c/Tarfia, sn \\ 
41012-Sevilla \\ 
Spain.} 




\email{arias@us.es, ariasdereyna1947@gmail.com}


\begin{abstract}
To define an explicit regions without zeros of $\Rzeta(s)$,  in a previous paper we obtained an approximation to $\Rzeta(s)$ of type $f(s)(1+U)$ with $|U|< 1$. But this $U$ do not tend to zero when $t\to+\infty$. In the present paper we get an approximation of the form $f(s)(1+\orden(t))$. We precise here Siegel's result, following his reasoning. This is essential to get the last Theorems in Siegel's paper about $\Rzeta(s)$.
\end{abstract}

\maketitle

\section{Introduction}

Siegel \cite{Siegel} studied the function 
\begin{equation}
\Rzeta(s)=\int_{0\swarrow1}\frac{x^{-s} e^{\pi i x^2}}{e^{\pi i x}-
e^{-\pi i x}}\,dx,
\end{equation}
that he found in Riemann's posthumous papers. Siegel related the zeros of $\Rzeta(s)$ to those of Riemann's zeta function.

In two previous papers \cite{A98}, \cite{A100}  we extend and  precise the asymptotic expansion of the auxiliary function of Riemann given in Siegel \cite{Siegel}. In \cite{A98} in order to get explicit regions without zeros of $\Rzeta(s)$, we obtain result of type 
\begin{equation}\label{defU2}
\Rzeta(s)=-\chi(s)\eta^{s-1}e^{-\pi i \eta^2}\frac{\sqrt{2}e^{3\pi i/8}\sin\pi\eta}
{2\cos2\pi\eta}(1+U),\qquad |U|<1,
\end{equation}
But these $U$ are not $\orden(t)$. To prove the last theorems in \cite{Siegel} about the zeros of $\Rzeta(s)$ and its connection to the zeros of $\zeta(s)$,  we need this type of approximation (see \cite{A102}). There are many possible results of this type valid in different regions $\Omega$, larger regions correspond to larger $U$.

We follow the steps of Siegel in this paper. Our main objective is to prove the  bounds $|U|\le A t^{-\beta}$ with absolute constants $A$ and $\beta>0$, and to stress the fact that $A$ is independent of $\Re(s)$ for regions of type $1-\sigma\ge t^\varepsilon$.

Commenting in the~\texttt{tex} file the line \texttt{excludecomment$\{A\}$} some proofs appear at the end of the paper, after the references, which we consider too elementary for inclusion on the paper.

\section{Asymptotic expansion}

\begin{proposition}\label{P:uno}
For $s=\sigma+it$ with $t>0$, $\sigma<0$ and $|s-1|>4\pi$ we have 
\begin{multline}\label{main}
\Rzeta(s)=\chi(s)\Bigl\{\sum_{n=k+1}^{\infty}n^{s-1}-\eta^{s-1}e^{-\pi i \eta^2}
\Bigl[\frac{\sqrt{2}e^{3\pi i/8}\sin\pi\eta-(-1)^k e^{2\pi i\eta-2\pi i(\eta-k)^2}}
{2\cos2\pi\eta}+\\+\sum_{j=m+1}^k(-1)^{j-1}e^{-2\pi i(j-\eta)^2}
w(j-\eta)+R\Bigr]\Bigr\},
\end{multline}
where 
$\eta=\sqrt{\frac{s-1}{2\pi i}}, \quad \eta_1=\Re(\eta),\quad \eta_2=\Im(\eta),\quad \eta_1+\eta_2>0$, and 
\[|R|\le\frac{\min(4, 29e^{-\pi \eta_2})}{|\eta|}+ 15\frac{e^{-\frac{\pi}{32}|\eta|^{2}}}{|\eta|},\]
$k$ is any integer $\ge m:=\lfloor\eta_1+\eta_2\rfloor$, (taking the sum in $j$ as $0$ in the case $k=m$).
Also, $w(z)$ defined in \cite{A98}*{eq.~(3)} satisfies the inequality 
\begin{equation}\label{E:boundw}
|w(x-\eta)|\le e^{\frac{4\pi}{3}\frac{|x-\eta|^3}{|\eta|}}-1, \quad \text{for $|x-\eta|\le \frac12|\eta|$}.
\end{equation}
\end{proposition}
\begin{proof}
The equation \eqref{main} is proved in  \cite{A98}*{eq.~(4)}, since $\sigma<0$ the zeta function there can be substituted by the Dirichlet series.

The bound of $R$ is proved \cite{A98}*{Thm.~5}.  The bound for $w(z)$ is \cite{A98}*{eq.~(23)} that is applied there to prove Theorem 5. 
\end{proof}

Let us define a closed set $\Omega\subset\C$ which contains the regions where our main Theorems will be true.  So our results in $\Omega$ applies for all our Theorems. 
Let  $0<\varepsilon <\frac12$ be given. 
\begin{equation}
\Omega:=\{s\in\C\colon\quad s=\sigma+it,\quad t\ge1,\quad 1-\sigma\ge t^\varepsilon,\quad  |s-1|\ge 4\pi\}.
\end{equation}

\begin{proposition}\label{P:eta}
Let $s$ be a point in $\Omega$ and $\eta=\eta_1+i\eta_2$ the point defined in the proposition \ref{P:uno}. Then we have $\sigma<0$,  $0<\eta_2< \eta_1$ and 
\begin{equation}\label{E:eta1}
\eta_1^2-\eta_2^2=\frac{t}{2\pi},\quad 2\eta_1\eta_2=\frac{1-\sigma}{2\pi}, \end{equation}
\begin{equation}\label{E:eta2}
\eta_1^2=\frac{\sqrt{t^2+(1-\sigma)^2}+t}{4\pi},\quad 
\eta_2^2=\frac{\sqrt{t^2+(1-\sigma)^2}-t}{4\pi}.
\end{equation}
\end{proposition}
\begin{proof}
For $s\in\Omega$ we have $t\ge1$ and $1-\sigma\ge 1$. Hence $\frac{s-1}{2\pi i}=\frac{t}{2\pi}+i\frac{1-\sigma}{2\pi}$ have argument between $0$ and $\pi/2$. Therefore, it has a square root with an argument between $0$ and $\pi/4$. This square root $\eta_1+i\eta_2$, satisfies the condition $\eta_1+\eta_2>0$, which defines $\eta$. It follows that $0<\eta_2<\eta_1$ (with strict inequalities because $(s-1)/(2\pi i)$ is not real or purely imaginary). 

Equation \eqref{E:eta1} follows by equating real and imaginary parts of $(\eta_1+i\eta_2)^2=(s-1)/(2\pi i)$. \eqref{E:eta2} follows from this easily.

Finally $1-\sigma\ge t^\varepsilon$ with $t\ge1$ implies $\sigma<0$ except for 
the case $t=1$ and $\sigma=0$. But the point $s=i$ is not contained in $\Omega$. 
\end{proof}

\begin{proposition}\label{P:tres}
For $s\in\Omega$ we have 
\begin{equation}\label{factor1}
\Bigl|-\frac{2\cos2\pi\eta}{\sqrt{2}e^{3\pi i/8}\sin\pi\eta}\Bigr|\le 
\frac{\sqrt{2}(1+\pi\eta_2)}{\pi\eta_2}e^{\pi\eta_2}.
\end{equation}
\end{proposition}
\begin{proof}
It is easy to see that $|\cos2\pi \eta|\le \cosh2\pi\eta_2$ and $|\sin\pi\eta|\ge\sinh\pi\eta_2$.

Since $\frac{x\cosh2x}{\sinh x}<(1+x)e^x$ for $x>0$,  we have
\begin{equation}
\Bigl|-\frac{2\cos2\pi\eta}{\sqrt{2}e^{3\pi i/8}\sin\pi\eta}\Bigr|\le 
\frac{2\cosh2\pi\eta_2}{\sqrt{2}\sinh\pi\eta_2}<\frac{\sqrt{2}}{\pi\eta_2}
(1+\pi\eta_2)e^{\pi\eta_2}.\qedhere
\end{equation}
\end{proof}

\begin{proposition}\label{P:cuatro}
For $\eta\in \Omega$ we have 
\begin{equation}
|\eta^{1-s}e^{\pi i \eta^2}|\le |\eta|^{1-\sigma}.
\end{equation}
\end{proposition}
\begin{proof}
By Proposition \ref{P:eta},  $0<\arg(\eta)<\frac{\pi}{4}$ we use this main argument to define $\eta^{1-s}$. Then we have
\[|\eta^{1-s}e^{\pi i \eta^2}|=\exp\bigl((1-\sigma)\log|\eta|+t\arg(\eta)-2\pi\eta_1\eta_2\bigr).\]
And we have 
\[t\arg(\eta)-2\pi\eta_1\eta_2=t\arctan\frac{\eta_2}{\eta_1}-2\pi\eta_1\eta_2
\le 2\pi\frac{t}{2\pi}\frac{\eta_2}{\eta_1}-2\pi\eta_1\eta_2.\]
Hence, by \eqref{E:eta1} 
\[t\arg(\eta)-2\pi\eta_1\eta_2\le 2\pi(\eta_1^2-\eta_2^2)\frac{\eta_2}{\eta_1}
-2\pi\eta_1\eta_2=-2\pi\frac{\eta_2^3}{\eta_1}<0.\]
And the result follows.
\end{proof}

\begin{proposition}\label{P:interm}
Let $s$ be a point in  $\Omega$, $m=\lfloor\eta_1+\eta_2\rfloor$ and $r\ge0$ an integer, then 
\begin{equation}\label{E:withU}
\Rzeta(s)=-\chi(s)\eta^{s-1}e^{-\pi i \eta^2}\frac{\sqrt{2}e^{3\pi i/8}\sin\pi\eta}{2\cos2\pi\eta}(1+U),
\end{equation}
where $U=U_1+U_2+U_3+U_4$ defined by 
\begin{gather}
U_1=-\frac{(-1)^{m+r}e^{2\pi i\eta-2\pi i (\eta-m-r)^2}}{\sqrt{2}e^{3\pi i/8}\sin\pi \eta},\qquad 
U_2=-\frac{2\cos2\pi\eta}{\sqrt{2}e^{3\pi i/8}\sin\pi \eta}\eta^{1-s}e^{\pi i\eta^2}\sum_{n>m+r}\frac{1}{n^{1-\sigma-it}},\\
U_3=\frac{2\cos2\pi\eta}{\sqrt{2}e^{3\pi i/8}\sin\pi \eta}
\sum_{j=m+1}^{m+r}(-1)^{j-1}e^{-2\pi i(j-\eta)^2}w(j-\eta),\qquad 
U_4=\frac{2\cos2\pi\eta}{\sqrt{2}e^{3\pi i/8}\sin\pi \eta} R.
\end{gather}
\end{proposition}
\begin{proof}
Equation \eqref{main} can be written as 
\begin{multline*}\Rzeta(s)=-\chi(s)\eta^{s-1}e^{-\pi i \eta^2}\frac{\sqrt{2}e^{3\pi i/8}\sin\pi\eta}{2\cos2\pi\eta}\Bigl\{1-\frac{(-1)^ke^{2\pi i\eta-2\pi i (\eta-k)^2}}{\sqrt{2}e^{3\pi i/8}\sin\pi \eta}\\+\frac{2\cos2\pi\eta}{\sqrt{2}e^{3\pi i/8}\sin\pi \eta}\Bigl[-\eta^{1-s}e^{\pi i\eta^2}
\sum_{n=k+1}^\infty\frac{1}{n^{1-\sigma-it}}+\sum_{j=m+1}^k(-1)^{j-1}e^{-2\pi i(j-\eta)^2}w(j-\eta)+R\Bigr]\Bigr\}
\end{multline*}
We put $k=m+r$ and obtain our result.
\end{proof}

\section{The approximation for \texorpdfstring{$1-\sigma$}{1-s} large}

\begin{theorem}\label{T:thfirst}
Let $\delta$ be a real number $\delta>\frac12$, and let $\Omega_1=\{s\in\Omega\colon 1-\sigma\ge t^{\delta}\ge2\}$. Then there is an absolute constant
$A$ such that for $s\in \Omega_1$
\begin{equation}\label{E:withUT1}
\Rzeta(s)=-\chi(s)\eta^{s-1}e^{-\pi i \eta^2}\frac{\sqrt{2}e^{3\pi i/8}\sin\pi\eta}{2\cos2\pi\eta}(1+U),\qquad |U|\le At^{-1/2}.
\end{equation}
\end{theorem}

\begin{proof}
We have \eqref{E:withU}, we will take $r=0$ so that $U_3=0$ and it remains to bound $U_1$, $U_2$ and $U_4$. 

For a point $s\in\Omega_1$ we have $1-\sigma\ge t^{\delta}$ and then by \eqref{E:eta1} and \eqref{E:eta2}
\begin{equation}\label{E:boundeta1}
\sqrt{\frac{t}{2\pi}}\le\eta_1\le |\eta|\le \Bigl(\frac{1-\sigma+t}{2\pi}\Bigr)^{1/2},
\end{equation}
and
\begin{equation}\label{E:boundeta2}
\eta_2=\frac{1-\sigma}{4\pi\eta_1}\ge \frac{1-\sigma}{2\sqrt{2\pi}(1-\sigma+t)^{1/2}}\ge\frac{t^\delta}{2\sqrt{2\pi}(t^\delta+t)^{1/2}}
\ge\frac{t^{\delta-\frac12}}{4\sqrt{\pi}} \gg t^{\delta-\frac12}.\end{equation}

By the definition of $m$ we have $m\le \eta_1+\eta_2<m+1$ so that
\begin{equation}\label{lastterm1}
|(-1)^{m-1}e^{2\pi i\eta-2\pi i(\eta-m)^2}|=e^{-2\pi\eta_2-4\pi(m-\eta_1)\eta_2}<e^{-2\pi\eta_2-4\pi(\eta_2-1)\eta_2}=e^{-2\pi(2\eta_2^2-\eta_2)}.
\end{equation}
Therefore,
\begin{equation}\label{E:U1bound1}
|U_1|=\Bigl|\frac{(-1)^ke^{2\pi i\eta-2\pi i (\eta-m-r)^2}}{\sqrt{2}e^{3\pi i/8}\sin\pi \eta}\Bigr|\le \frac{e^{-2\pi(2\eta_2^2-\eta_2)}}{\sqrt{2}\sinh\pi\eta_2}\le  
\frac{4}{\sqrt{2\pi}}t^{\frac12-\delta} \exp(-\tfrac14 t^{2\delta-1}+\tfrac{\sqrt{\pi}}{2}t^{\delta-\frac12}).\end{equation}
To bound $U_2$, we apply the propositions \ref{P:tres} and \ref{P:cuatro}, and then bound the sum in absolute value by the first term and an integral
\[|U_2|\le \frac{\sqrt{2}(1+\pi\eta_2)}{\pi\eta_2}e^{\pi\eta_2}|\eta|^{1-\sigma}
\Bigl\{\frac{(m+1)^\sigma}{-\sigma}+(m+1)^{\sigma-1}\Bigr\}.\]

We have $m=\lfloor \eta_1+\eta_2\rfloor$ and $\sigma<0$, therefore (by applying 
\eqref{E:eta1})
\begin{align*}
&\Bigl(\frac{|\eta|}{m+1}\Bigr)^{1-\sigma}\Bigl(1+\frac{m+1}{-\sigma}\Bigr)\le \Bigl(\frac{|\eta|}{\eta_1+\eta_2}\Bigr)^{1-\sigma}\Bigl(1+\frac{\eta_1+\eta_2+1}{-\sigma}\Bigr)\\
&= \Bigl(\frac{\eta_1^2+\eta_2^2}{\eta_1^2+\eta_2^2+2\eta_1\eta_2}\Bigr)^{\frac{1-\sigma}{2}}\Bigl(1+\frac{\eta_1+\eta_2+1}{4\pi\eta_1\eta_2-1}\Bigr)\\
&=\Bigl(1+\frac{\eta_1+\eta_2+1}{4\pi\eta_1\eta_2-1}\Bigr)\exp\Bigl\{-2\pi\eta_1\eta_2\log\Bigl(1+\frac{2\eta_1\eta_2}{\eta_1^2+\eta_2^2}\Bigr)\Bigr\}.
\end{align*}
For  $0<x=\frac{2\eta_1\eta_2}{\eta_1^2+\eta_2^2}\le 1$ we have $\log(1+x)\ge x/2$, so that 
\[
|\eta|^{1-\sigma}
\Bigl\{\frac{(m+1)^\sigma}{-\sigma}+(m+1)^{\sigma-1}\Bigr\}\le \Bigl(1+\frac{\eta_1+\eta_2+1}{4\pi\eta_1\eta_2-1}\Bigr)e^{-2\pi\frac{\eta_1^2\eta_2^2}{\eta_1^2+\eta_2^2}}.
\]
We have $1\le t$ and by \eqref{E:boundeta1} and \eqref{E:boundeta2}
\[\Bigl(1+\frac{\eta_1+\eta_2+1}{4\pi\eta_1\eta_2-1}\Bigr)\le 1+\frac{\eta_1+\eta_2+1}{2\pi\eta_1\eta_2}\le1+\frac{4\sqrt{\pi}}{2\pi} +\frac{\sqrt{2\pi}}{2\pi}+\frac{4\pi\sqrt{2}}{2\pi}=a<\frac{11}{2}.\]
The quantity in parentheses is bounded by a constant and since $\eta_2<\eta_1$ we obtain 
\begin{equation}\label{E:boundU21}
|U_2|\le \sqrt{2}(1+\tfrac{1}{\pi\eta_2})\;a e^{\pi \eta_2-\pi \eta_2^2}\le 8(1+\tfrac{1}{\pi\eta_2})\exp\bigl(-\tfrac{1}{16} t^{2\delta-1}+\tfrac{\sqrt{\pi}}{4}t^{\delta-\frac12}\bigr)\ll \exp\bigl(-\tfrac{1}{32}t^{2\delta-1}).\end{equation}
By Propositions \ref{P:uno} and \eqref{factor1} we have 
\begin{equation}\label{E:boundU41}
\begin{aligned}
|U_4|&\le \frac{\sqrt{2}(1+\pi\eta_2)}{\pi\eta_2}e^{\pi\eta_2}\Bigl(\frac{29e^{-\pi \eta_2}}{|\eta|}+ 15\frac{e^{-\frac{\pi}{32}|\eta|^{2}}}{|\eta|}\Bigr)\\
&\le\sqrt{2}\Bigl(1+\frac{1}{\pi\eta_2}\Bigr)\Bigl(\frac{29}{|\eta|}+\frac{15e^{\pi\eta_2-\frac{\pi}{32}\eta_2^2}e^{-\frac{\pi}{32}\eta_1^2}}{|\eta|}\Bigr)\le \frac{C}{|\eta|}
\le \frac{C}{\eta_1}\ll t^{-1/2}.
\end{aligned}
\end{equation}
The function $e^{\pi\eta_2-\frac{\pi}{32}\eta_2^2}$ is bounded by $e^{8\pi}$, but decreases exponentially for $\eta_2\to+\infty$. For $s\in\Omega_2$ we have 
$\eta_2\gg t^{\delta-\frac12}\to+\infty$.
Since $U_1$, and $U_2$ decrease exponentially for $t\to+\infty$, $U$ is of the order of $U_4$.
\end{proof}

\begin{remark}
Taking, for example $\delta=\frac59$ in Theorem \ref{T:thfirst}, it is true that $U\le At^{-1/2}$, but the constant $A$ is very large. The limit $\lim_{t\to+\infty}t^{1/2}U=\lim_{t\to+\infty}t^{1/2}U_4=58\sqrt{\pi}$. But the size of $A$ is determined by the other terms; for example, $t^{1/2}U_1$ increases until $t\approx 10^{13}$ where $t^{1/2}U_1\approx 97697$. 
\end{remark}

\section{Approximation with \texorpdfstring{$1-\sigma$}{1-sigma} moderate}

In this section $t^\delta\ge1-\sigma\ge t^\varepsilon$ with $ \varepsilon<\frac12<\delta$. For a given $t$, no matter how large, the value of $\eta_2=(1-\sigma)/4\pi\eta_1$ can be $\eta_2=16$, so that $e^{\pi \eta_2-\pi\eta_2^2/32}=e^{8\pi}$ is large or even very small. This has to be considered in what follows.

\begin{theorem}\label{T:second}
Given  $\varepsilon$ with $\frac25<\varepsilon<\frac12$, then there are real numbers $\delta>0$ and $\alpha>0$ such that 
\[\tfrac25<\varepsilon<\tfrac12<\delta<1-\varepsilon<1,\quad 
0<\tfrac12-\varepsilon<\alpha<\tfrac\varepsilon4<\tfrac16.\]
With $\varepsilon$, $\delta$ and $\alpha$ fixed as above, there is a $t_0\ge 1$ and an absolute constant $A$ such that for $s=\sigma+it$ with $t\ge t_0$ and $t^\delta\ge1-\sigma\ge t^{\varepsilon}$ we have 
\begin{equation}\label{E:withUT2}
\Rzeta(s)=-\chi(s)\eta^{s-1}e^{-\pi i \eta^2}\frac{\sqrt{2}e^{3\pi i/8}\sin\pi\eta}{2\cos2\pi\eta}(1+U),\qquad |U|\le At^{4\alpha-\varepsilon}.
\end{equation}
\end{theorem}
\begin{proof}
Let $\Omega_2:=\{s\in\Omega\colon t^{\delta}\ge1-\sigma\ge t^{\varepsilon},t\ge t_0\}$.
In what follows, we assume that $s\in\Omega_2$. In this case $\eta$ is given by a convergent power series
\begin{equation}\label{etat}
\eta=\Bigl(\frac{t}{2\pi}+\frac{1-\sigma}{2\pi}i\Bigr)^{1/2}=
\sqrt{\frac{t}{2\pi}}\Bigl(1+i\frac{1-\sigma}{2t}+\frac{(1-\sigma)^2}{8t^2}-i\frac{(1-\sigma)^3}{16t^3}-\frac{5(1-\sigma)^4}{128t^4}+\cdots\Bigr)
\end{equation}
so that
\begin{equation}\label{eta12t}
\eta=\Bigl(\frac{t}{2\pi}\Bigr)^{1/2}+\Orden\Bigl(\frac{\sigma-1}{t^{1/2}}\Bigr),\quad
\eta_1=\Bigl(\frac{t}{2\pi}\Bigr)^{1/2}+\Orden\Bigl(\frac{(1-\sigma)^2}{t^{3/2}}\Bigr),\quad \eta_2=\frac{1-\sigma}{2\sqrt{2\pi t}}+\Orden\Bigl(\frac{(1-\sigma)^3}{t^{5/2}}\Bigr).
\end{equation}
By \eqref{E:eta1} and  \eqref{E:eta2} 
\[\frac{t}{2\pi}\le \eta_1^2\le\frac{\sqrt{t^2+t^{2\delta}}+t}{4\pi}\le \frac{\sqrt{2}+1}{4\pi}t,\quad\text{and}\quad \eta_2=\frac{1-\sigma}{4\pi\eta_1}\le  \frac{t^{\delta-\frac12}}{2\sqrt{2\pi}},\]
it follows that for $s\in\Omega_2$,
\begin{equation}\label{E:etasumbound}
\eta_1+\eta_2\le a_1 t^{1/2}, \quad a_1:=\frac{\sqrt{\sqrt{2}+1}+2^{-1/2}}{2\sqrt{\pi}}.
\end{equation}

For $s$ in $\Omega_2$ we will take, in Proposition \ref{P:interm} the integer $r\ge1$ as a function of $t$. Our objective is to bound the $U_i$.
In the course of our reasoning, we will impose several conditions on $r$, which at the end will be satisfied with an adequate election of $r$.

Since $m\le\eta_1+\eta_2<m+1$ we have 
\begin{displaymath}
|e^{2\pi i\eta-2\pi i(m+r-\eta)^2}|=e^{-2\pi \eta_2-4\pi(m+r-\eta_1)\eta_2}<e^{-2\pi \eta_2-4\pi (r+\eta_2-1)\eta_2}<e^{-2\pi r \eta_2},
\end{displaymath}
we have used here the condition $r\ge1$ to get a term $-2\pi\eta_2$ in the exponent. So our first condition on $r$
\begin{equation}\tag{{$\bullet$}}
r\ge 1.
\end{equation}
We also have  
\begin{equation}\label{E:interm}
|\sin\pi\eta|\ge\sinh\pi\eta_2\ge\pi\eta_2\ge \frac{\pi}{4\pi\frac{\sqrt{\sqrt{2}+1}}{2\sqrt{\pi}}}\frac{1-\sigma}{\sqrt{t}}\ge \frac{\sqrt{\pi}}{2\sqrt{\sqrt{2}+1}}t^{\varepsilon-\frac12}= a_2 t^{\varepsilon-\frac12}.
\end{equation}
So, 
\[|U_1|=\Bigl|-\frac{(-1)^{m+r}e^{2\pi i\eta-2\pi i (\eta-m-r)^2}}{\sqrt{2}e^{3\pi i/8}\sin\pi \eta}\Bigr|\le\frac{1}{a_2\sqrt{2}} \; t^{\frac12-\varepsilon}e^{-2\pi r\eta_2}.\]
To bound $U_2$, we proceed as in the proof of Theorem \ref{T:thfirst} applying Propositions \ref{P:tres} and \ref{P:cuatro}, and then bound the sum in absolute value by the first term and an integral
\begin{multline}\label{E:U2bound}
|U_2|=\Bigl|\frac{2\cos2\pi\eta}{\sqrt{2}e^{3\pi i/8}\sin\pi \eta}\eta^{1-s}e^{\pi i\eta^2}\sum_{n>m+r}\frac{1}{n^{1-\sigma-it}}\Bigr|\\
\le
\frac{\sqrt{2}(1+\pi\eta_2)}{\pi\eta_2}e^{\pi\eta_2}|\eta|^{1-\sigma}
\Bigl\{\frac{(m+r+1)^\sigma}{-\sigma}+(m+r+1)^{\sigma-1}\Bigr\}\\=
\frac{\sqrt{2}(1+\pi\eta_2)}{\pi\eta_2}e^{\pi\eta_2}\Bigl(1+\frac{m+r+1}{-\sigma}\Bigr)\Bigl(\frac{|\eta|}{m+r+1}\Bigr)^{1-\sigma}.
\end{multline}

Note that 
\[|\eta|^2=\Bigl|\frac{s-1}{2\pi i}\Bigr|=\frac{((1-\sigma)^2+t^2)^{1/2}}{2\pi}\le \frac{(t^{2\delta}+t^2)^{1/2}}{2\pi}\le \frac{(2t^2)^{1/2}}{2\pi}=\sqrt{2}\frac{t}{2\pi}\le \sqrt{2}\eta_1^2,\]
and therefore $|\eta|\le \frac43\eta_1$, (check $2^{1/4}\le 4/3$). In addition, to continue the bound of $|U_2|$ we need another hypothesis on $r$
\begin{equation}\tag{${\bullet}{\bullet}$}
r\le |\eta|.
\end{equation}
Since $m\le \eta_1+\eta_2<m+1$ and $|\eta|\le \eta_1+\eta_2<m+1$, we have $\frac{|\eta|+r}{|\eta|}<\frac{m+r+1}{|\eta|}$ and then 
(~notice that $1+x\ge e^{2x/3}$ for $0\le x\le1$~)
\[\Bigl(\frac{m+r+1}{|\eta|}\Bigr)^{\sigma-1}<\Bigl(\frac{|\eta|+r}{|\eta|}\Bigr)^{\sigma-1}\le e^{(\sigma-1)\frac{2r}{3|\eta|}}=e^{-(1-\sigma)\frac{2r}{3|\eta|}}<e^{-(1-\sigma)\frac{r}{2\eta_1}}=e^{-2\pi r\eta_2}.\]
And by \eqref{E:interm} and \eqref{E:U2bound} it follows that 
\begin{equation}
\begin{aligned}
|U_2|&\le \sqrt{2}(1+\tfrac{1}{\pi\eta_2})\Bigl(1+\frac{\eta_1+\eta_2+r+1}{-\sigma}\Bigr)e^{\pi\eta_2-2\pi r\eta_2}\\&\le \sqrt{2}\bigl(1+\tfrac{1}{a_2}t^{\frac12-\varepsilon}\bigr)\Bigl(1+\frac{a_1t^{1/2}+r+1}{t^\varepsilon-1}\Bigr)
e^{-\pi r\eta_2}.
\end{aligned}
\end{equation}
We start now the bound of $U_3$ and this will be the largest term. To bound $w(j-\eta)$ we have to apply \eqref{E:boundw}. Hence, for $1\le \ell\le r$ we need 
$|m+\ell-\eta|\le \frac12|\eta|$. Since
\begin{displaymath}
|m+\ell-\eta|^2=(m+\ell-\eta_1)^2+\eta_2^2\le (r+\eta_2)^2+\eta_2^2
\le 2r^2+3\eta_2^2, 
\end{displaymath}
we add a condition to $r$
\begin{equation}\tag{{$\bullet$}{$\bullet$}{$\bullet$}}
2r^2+3\eta_2^2\le \frac14|\eta|^2.
\end{equation}
Then, assuming 
\begin{equation}\tag{${\bullet}{\bullet}{\bullet}{\bullet}$}
\eta_2\le r, 
\end{equation}
we get 
\begin{equation}
|w(m+\ell-\eta)|\le e^{\frac{4\pi}{3}\frac{(5r^2)^{3/2}}{|\eta|}}-1\le \frac{8\pi 5^{3/2}}{3}\frac{r^3}{|\eta|},
\end{equation}
because $e^x-1\le 2x$ for $0\le x\le1$ and we assume 
\begin{equation}\tag{${\bullet}{\bullet}{\bullet}{\bullet}{\bullet}$}
\frac{4\pi(5^{3/2}) r^3}{3|\eta|}\le 1.
\end{equation}
Joining this with \eqref{factor1} and 
\begin{displaymath}
|e^{-2\pi i(m+\ell-\eta)^2}|\le e^{-4\pi(\eta_2+\ell)\eta_2}\le e^{-4\pi\eta_2^2}
\end{displaymath}
yields
\[|U_3|\le \sqrt{2}\bigl(1+\tfrac{1}{a_2}t^{\frac12-\varepsilon}\bigr)r e^{\pi\eta_2-4\pi\eta_2^2}\frac{8\pi 5^{3/2}}{3}\frac{r^3}{|\eta|}.\]
Notice that for $s\in \Omega_2$ the value of $\eta_2\sim (1-\sigma) t^{-1/2}$ range from very small to large values. So that $e^{\pi\eta_2-4\pi\eta_2^2}\le 5/4$ is bounded only by the value when $\eta_2=1/8$. Therefore, 
\begin{equation}
|U_3|\le a_3\bigl(1+\tfrac{1}{a_2}t^{\frac12-\varepsilon}\bigr)\frac{r^4}{|\eta|},\qquad a_3=\tfrac13 40\pi\sqrt{10}e^{\pi/16}.
\end{equation}
The bound of $U_4$ is similar to that in Theorem \ref{T:thfirst}
\begin{align*}
|U_4|&\le \sqrt{2}\Bigl(1+\frac{1}{\pi\eta_2}\Bigr)\Bigl(\frac{29}{|\eta|}+\frac{15e^{\pi\eta_2-\frac{\pi}{32}\eta_2^2-\frac{\pi}{32}\eta_1^2}}{|\eta|}\Bigr)\\
&\le \sqrt{2}\bigl(1+\tfrac{1}{a_2}t^{\frac12-\varepsilon}\bigr)\bigl(29+15e^{8\pi-\frac{t}{64}}\bigr)\sqrt{2\pi}\;t^{-1/2}.
\end{align*}
For $s\in \Omega_2$ there is no better bound than $e^{\pi\eta_2-\frac{\pi}{32}\eta_2^2}\le e^{8\pi}$, even when $t$ is large since $\eta_2\asymp (1-\sigma)t^{-1/2}$. 

In summary, for $t\to+\infty$ and with absolute constants 
\[|U|\ll t^{\frac12-\varepsilon}e^{-2\pi r\eta_2}+t^{\frac12-\varepsilon}(t^{\frac12-\varepsilon}+r t^{-\varepsilon})e^{-\pi r \eta_2}+\frac{r^4}{t^{\varepsilon}}
+t^{-\varepsilon}.\]
The first two terms will be controlled if  $r\eta_2\gg rt^{\varepsilon-\frac12}$ tend to infinity for $t\to+\infty$ and we also need the term $r^4t^{-\varepsilon}\to0$ for $t\to+\infty$. This will be only possible if 
for $t$ large we have 
\[rt^{\varepsilon-\frac12}>rt^{-\frac\varepsilon4}.\]
Hence, for $\varepsilon >2/5$. Assuming this, we take $r=\lfloor t^\alpha\rfloor$
with 
\[0<\tfrac12-\varepsilon<\alpha<\tfrac\varepsilon4.\]
Taking any of these values of $\alpha$, there is $t_0\ge1$ such that all the conditions on $r$  are true. The first is true since for $t\ge1$, the value of 
$r=\lfloor t^\alpha\rfloor\ge1$. For the second since $\alpha<\frac18$ we have 
\[r\le t^\alpha\le (t/2\pi)^{1/2}\le |\eta|, \qquad t\ge t_0.\]
The third condition follows from 
\[2r^2+3\eta_2^2\le 2t^{2\alpha}+3\frac{(1-\sigma)^2}{16\pi^2\eta_1^2}\ll 
t^{2\alpha}+t^{2\delta-1}=\orden(t), \]
while $|\eta|^2\gg t$. 
The fourth condition follows from $\eta_2\ll t^{\delta-1/2}$ and $r=t^\alpha$
and $\delta-\frac12<\frac12-\varepsilon<\alpha$. 
The fifth condition follows from 
\[\frac{4\pi r^3}{|\eta|}\ll t^{3\alpha-\frac12}\to0.\]

With these choices, we have for any $s\in\Omega_2$ and with absolute constants
\[|U|\ll t^{\frac12-\varepsilon}e^{-c_1 t^{\alpha+\varepsilon-\frac12}}+t^{\frac12-\varepsilon}(t^{\frac12-\varepsilon}+r t^{-\varepsilon})e^{-c_2 t^{\alpha+\varepsilon-\frac12}} +t^{4\alpha-\varepsilon}
+t^{-\varepsilon}\ll t^{4\alpha-\varepsilon}.\qedhere\]
\end{proof}

\begin{theorem}
There exist constants $A$ and $t_0>1$ such that for $s$ in the closed set
\[\Omega=\{s\in\C\colon t\ge t_0,\quad 1-\sigma\ge t^{3/7}\},\]
we have 
\begin{equation}\label{E:withUT3}
\Rzeta(s)=-\chi(s)\eta^{s-1}e^{-\pi i \eta^2}\frac{\sqrt{2}e^{3\pi i/8}\sin\pi\eta}{2\cos2\pi\eta}(1+U),\qquad |U|\le At^{-\frac{1}{21}}.
\end{equation}
\end{theorem}

\begin{proof}
In Theorem \ref{T:second} take $\varepsilon=3/7$, with $\frac25<\frac37<\frac12$.
Then we may pick $\delta=\frac59$ with $\frac12<\frac59<\frac47$, and $\alpha=\frac{2}{21}$ satisfying $\frac12-\varepsilon=\frac{1}{14}<\frac{2}{21}<\frac{3}{28}=\frac\varepsilon4$. We apply Theorem \ref{T:thfirst} when $1-\sigma\ge t^\delta$, obtaining an error $U$ bounded by $At^{-1/2}$ and Theorem 
\ref{T:second} for $t^\delta\ge1-\sigma\ge t^{\varepsilon}$, obtaining an error 
bounded by $A' t^{-1/21}$, where $-\frac{1}{21}=4\alpha-\varepsilon$.
\end{proof}

\begin{A}
\begin{proof}[Proof of the inequality $\frac{x\cosh2x}{\sinh x}<(1+x)e^x$]
The following inequalities are all equivalent and the last is trivial  for $x>0$. \begin{gather*}
x(e^{2x}+e^{-2x})<(1+x)e^x(e^x-e^{-x})\\
x(e^{2x}+e^{-2x})<(1+x)(e^{2x}-1)\\
e^{2x}-1+xe^{2x}-x>xe^{2x}+xe^{-2x}\\
e^{2x}-1-x>xe^{-2x};\qquad
e^{2x}-1>x(1+e^{-2x})\\
x<e^{2x}\frac{1-e^{-2x}}{1+e^{-2x}}\quad \text{since}\quad
e^{2x}\frac{1-e^{-2x}}{1+e^{-x}}<e^{2x}\frac{1-e^{-2x}}{1+e^{-2x}}\\
\end{gather*}
it is sufficient to prove 
\[x<e^{2x}\frac{1-e^{-2x}}{1+e^{-x}}=e^{2x}(1-e^{-x})=e^{2x}-e^x.\]
But the inequality $x<e^{2x}-e^x$ is trivial (think of the power series).
\end{proof}
\end{A}


\begin{thebibliography}{999}
%
\bibitem{A98}
\textsc{J.~Arias de Reyna}, \emph{Region without zeros for the auxiliary function of Riemann}, \href{https://arxiv.org/abs/2406.03825}{arXiv:2406.03825}.  


\bibitem{A100}
\textsc{J.~Arias de Reyna}, \emph{Asymptotic expansions of the auxiliary function},  \href{https://arxiv.org/pdf/2406.04714}{arXiv:2406.04714}. 

\bibitem{A102}
\textsc{J.~Arias de Reyna}, \emph{On Siegel results about the zeros of the auxiliary function of Riemann},  preprint (102).

\bibitem{Siegel}
\textsc{C.~L.~Siegel}, \emph{Über Riemann Nachla\ss\ zur analytischen
Zahlentheorie}, Quellen und Studien zur Geschichte der Mathematik Astronomie und 
Physik \textbf{2} (1932) 45--80. (Reprinted in \cite{SW}, 1, 275--310.)
\href{https://arxiv.org/abs/1810.05198}{English version}.

\bibitem{SW}
\textsc{C. L. Siegel}, \emph{Carl Ludwig Siegel's Gesammelte Abhandlungen}, 
(edited by K. Chandrasekharan and H. Maa\ss), Springer-Verlag, Berlin, 1966.

\end{thebibliography}
\end{document}